\documentclass[12pt, reqno]{amsart}
\usepackage{ifthen,amsfonts,amsmath,amssymb,epic,eepic}
\usepackage{epsfig,color}

\makeatletter
\@addtoreset{equation}{section}
\makeatother

\theoremstyle{definition}

\newtheorem{Thm}{Theorem}
\newtheorem{Cor}{Corollary}

\newtheorem{Lem}{Lemma}

\newtheorem{Def}{Definition}

\newtheorem{Pro}{Proposition}

\renewcommand{\rm}{\normalshape} %makes rm mean roman

\newcommand{\nn}{{\bf{n}}}

\newcommand{\cS}{{\mathcal{S}}}
\newcommand{\cP}{{\mathcal{P}}}

\newcommand{\dist}{{\text {dist}}}

\newcommand{\K}{{\text K}}

\def\RR{{\bf  R}}

\def\CC{{\bf C }}
\def\cone{{\bf C }}

\newcommand{\Area}{{\text {Area}}}

\newcommand{\cB}{{\mathcal{B}}}

\newcommand{\cF}{{\mathcal{F}}}

\newcommand{\eqr}[1]{(\ref{#1})}

\begin{document}

\title[Disks that are double spiral staircases]
{Disks that are double spiral staircases}
%\title[Embedded minimal disks in $\RR^3$]
%{Embedded minimal disks in $\RR^3$}

\author{Tobias H. Colding}%
\address{Courant Institute of Mathematical Sciences and Princeton University\\
251 Mercer Street\\ New York, NY 10012 and Fine Hall, Washington
Rd., Princeton, NJ 08544-1000}
\author{William P. Minicozzi II}%
\address{Department of Mathematics\\
Johns Hopkins University\\
3400 N. Charles St.\\
Baltimore, MD 21218}
\thanks{The authors were partially supported by NSF Grants DMS
0104453 and DMS 0104187}

%%\date{\today}

\email{colding@cims.nyu.edu, minicozz@jhu.edu}

%\renewcommand{\abstractname}{Summary}
%\begin{abstract}
%\end{abstract}

\maketitle

%%%%%%%%%%%%%
%{\small
%\tableofcontents}
%%%%%%%%%%%%%

%\numberwithin{equation}{section}

%%numbering for equations

\setcounter{part}{0}
%\numberwithin{section}{part} %%%number sections within parts
\renewcommand{\rm}{\normalshape} %makes rm mean roman
\renewcommand{\thepart}{\Roman{part}}
\setcounter{section}{1}

\centerline{What are the possible shapes of various things and why?}
  
\vskip2mm
For instance, 
when a closed wire or a frame 
is dipped into a soap solution and is raised up from 
the solution, the surface spanning the wire is a soap film; 
see p. 11 of \cite{Op} or 
fig. I that show a soap film with the shape of a (single) spiral staircase.  
What are the 
possible shapes of soap films and why?  Or, for instance, 
why is DNA like a double spiral 
staircase?  ``What..?'' and ``why..?'' are fundamental questions, 
and when answered,  
help us understand the world we live in.  

Soap films, soap bubles, 
and surface tension were extensively studied by the Belgian physicist 
and inventor (the inventor of the stroboscope) Joseph Plateau 
in the first half of the 
nineteenth century.  At least since his studies, it has been known 
that the right  
mathematical model for soap films 
are minimal surfaces -- the soap film is in a state of 
minimum energy when it is covering the least possible amount of area.  

We will discuss here the answer to the question: ``What are the possible 
shapes of embedded minimal disks in $\RR^3$ and why?''.

The field of minimal surfaces  
dates back to the publication in 1762 of 
Lagrange's famous memoir ``Essai d'une nouvelle 
m\'ethode pour d\'eterminer les maxima et les minima des formules 
int\'egrales ind\'efinies''.  Euler had already in a paper published in 
1744 discussed minimizing properties of the surface now known as the 
catenoid, but he only considered variations within a 
certain class of surfaces.  
In the almost one quarter of a millenium 
that has past 
since Lagrange's memoir minimal 
surfaces has remained a vibrant area of research and 
there are many reasons why.  The study 
of minimal surfaces was the
birthplace of regularity theory.  It lies
on the intersection of nonlinear elliptic PDE, geometry, and
low--dimensional topology and over the years the field has matured 
through the efforts of many people.  However some very
fundamental questions remain.  Moreover, many of the potentially
spectacular applications of the field have yet to be achieved.  For instance, 
it has long been the hope that several of the outstanding conjectures 
about the topology of $3$-manifolds could be resolved using detailed 
knowledge of minimal surfaces.  Surfaces with
uniform curvature (or area) bounds have been well understood and the regularity
theory is complete, yet essentially nothing was known without such
bounds.  We discuss here the theory of embedded minimal disks in $\RR^3$ 
without a priori bounds.  As we will see, the helicoid, which is a double
spiral staircase, is the most important example of such a disk.  In 
fact, we will see that every embedded minimal disk is 
either a graph of a function or is part of a
double spiral staircase.  The helicoid was discovered to be a minimal surface 
by Meusnier in 1776.

\newpage
\vskip2mm
\noindent
{\bf{Double spiral staircases}}
\vskip2mm

A {\it double spiral staircase} consists of two staircases
  that spiral around one another so that without meeting two people can 
pass each other.  For 
instance, one could ascend one staircase while the other descends 
the other staircase.  Fig. II\\  
(http://www.a-castle-for-rent.com/castles/images/Chambord14.jpg)\\ 
shows Leonardo da Vinci's 
double spiral staircase
in Ch\^ateau de Chambord in the Loire valley in France.  
The construction of the 
castle began in 1519 (the same year that Leonardo da Vinci died) and was 
completed in 1539.  Fig. III\\
(http://www.angelfire.com/trek/lafrance1999/tours3.html)\\
show a model of the double spiral staircase where 
we can clearly see the two staircases spiraling around one another.  
In fig. IV\\ 
(http://images.amazon.com/images/P/0451627873.01.LZZZZZZZ.jpg)\\
we see ``the double helix'' which was 
discovered in 1952 by Crick and Watson to be the structure of DNA.  
The double spiral 
structure represented the culmination of half a century of prior 
work on genetics and is by many considered one of the 
greatest scientific discoveries of the twentieth century.  
Also the internal ear,
the cochlea, is a double spiral staircase; see fig. V or p. 343 of \cite{K}.

In the cochlea, the two canals wind around a conical bony axis and after
about two and a half rotations they meet at the top and fuse.  
The canals are filled with fluids and 
sound waves travel up one canal, turn around, and 
come down the other.  When the liquid is set into movement, it will set
the Basilar membrane and the hair cells into vibration.  Different hair cells
correspond to different frequencies.

Other examples of double spiral staircases include parking ramps.     

\vskip2mm
\noindent
{\bf{What is a minimal surface and what are the central examples?}}
\vskip2mm

Let $\Sigma\subset \RR^3$ be a \underline{smooth} 
orientable surface (possibly with
boundary) with unit normal $\nn_{\Sigma}$.  Given
a function $\phi$ in the space $C^{\infty}_0(\Sigma)$ 
of infinitely differentiable (i.e., smooth), compactly supported 
functions on $\Sigma$, consider the one-parameter variation
\begin{equation}
\Sigma_{t,\phi}=\{x+t\,\phi (x)\,\nn_{\Sigma}(x) | x\in \Sigma\}\, .
\end{equation}
The so called first variation formula of area is the equation (integration is
with respect to $d\text{area}$)
\begin{equation}  \label{e:frstvar}
\left.\frac{d}{dt} \right|_{t=0}\Area (\Sigma_{t,\phi})
=\int_{\Sigma}\phi\,H\, ,
\end{equation}
where $H$ is the mean curvature of $\Sigma$.  (When $\Sigma$ is noncompact, 
then $\Sigma_{t,\phi}$ in \eqr{e:frstvar} is replaced by 
$\Gamma_{t,\phi}$, where 
$\Gamma$ is any compact set containing the support of $\phi$.)  
The surface $\Sigma$ is said to be a {\it minimal} surface (or just minimal) if
\begin{equation}
\left.\frac{d}{dt} \right|_{t=0}\Area (\Sigma_{t,\phi})=0
\,\,\,\,\,\,\,\,\,\,\,\text{ for all } \phi\in C^{\infty}_0(\Sigma)
\end{equation}
 or, equivalently by \eqr{e:frstvar}, if the
mean curvature $H$ is identically 
zero.  Thus $\Sigma$ is minimal if and only if
it is a critical point for the area functional.  (Since a critical 
point is not necessarily a minimum the term 
``minimal'' is
misleading, but it is time honored.  The equation for a critical 
point is also sometimes called the Euler-Lagrange equation.)  
Moreover, a computation shows that if $\Sigma$ is minimal, then
\begin{equation}
\left. \frac{d^2}{dt^2} \right|_{t=0}\Area (\Sigma_{t,\phi})
=-\int_{\Sigma}\phi\,L_{\Sigma}\phi\, ,
\,\,\,\,\,\,\,\,\,\,\,\text{ where }
L_{\Sigma}\phi=\Delta_{\Sigma}\phi+|A|^2\phi
\end{equation}
is the second variational (or Jacobi) operator.
Here $\Delta_{\Sigma}$ is the Laplacian on $\Sigma$ and
$A$ is the second fundamental form.  So $|A|^2=\kappa_1^2+\kappa_2^2$,
 where $\kappa_1,\,\kappa_2$ are the principal curvatures of $\Sigma$
and $H=\kappa_1+\kappa_2$.
A minimal surface $\Sigma$ is said to be
stable if
\begin{equation}
\left. \frac{d^2}{dt^2} \right|_{t=0}\Area (\Sigma_{t,\phi})\geq 0
\,\,\,\,\,\,\,\,\,\,\,\text{ for all } \phi\in C^{\infty}_0(\Sigma)\, .
\end{equation}
One can show that a minimal graph is stable and, more generally, so is a
multi-valued minimal graph (see below for the precise definition).

There are two local models for embedded minimal disks (by an {\it
  embedded disk} we mean a smooth injective map from the closed unit ball in
$\RR^2$ into $\RR^3$).  One model is the plane
(or, more generally, a minimal graph) and the other is a piece of a 
helicoid.

The derivation of the equation for a minimal graph goes back 
to Lagrange's 1762 memoir.

\vskip2mm
\noindent
{\bf{Example 1}}: 
(Minimal graphs).  If $\Omega$ is a simply connected domain in $\RR^2$ and 
$u$ is a real valued function on $\Omega$ satisfying the minimal surface 
equation
\begin{equation}  \label{e:mineq}
\text{div} \left( \frac{\nabla u}{\sqrt{1+|\nabla u|^2}}\right)=0\, ,
\end{equation}
then the graph of $u$, i.e., the set $\{(x_1,x_2,u(x_1,x_2))\,|\,
(x_1,x_2)\in \Omega\}$, is a minimal disk.

\vskip2mm
A classical theorem of Bernstein from 1916 says that entire (i.e., where
$\Omega=\RR^2$) minimal graphs are planes.  This remarkable theorem of 
Bernstein was one of the first illustrations of the fact that the solutions 
to a nonlinear PDE, like the minimal surface equation, can behave 
quite differently from the solutions to a linear equation.   In the early 
nineteen--eigthties Schoen and Simon extended the theorem of 
Bernstein to complete 
simply connected embedded minimal surfaces in $\RR^3$ with quadratic area 
growth.  A surface $\Sigma$ is said to have quadratic area growth if  
for all $r>0$, the intersection of the surface with the ball in $\RR^3$ of 
radius $r$ and center at the origin is bounded by $C\, r^2$ for a 
fixed constant $C$ independent of $r$. 

The second model comes from the helicoid which was discovered by 
Meusnier in 1776.  Meusnier had been a student of Monge.  He also 
discovered that the surface now known as the catenoid is minimal 
in the sense of Lagrange, and he was  
the first to characterize a minimal surface as a surface with 
vanishing mean curvature.  Unlike the helicoid, the catenoid is 
not topologically a plane but rather a cylinder.  

The 
helicoid is a ``double spiral staircase''.

\begin{figure}[htbp]
    \setlength{\captionindent}{20pt}
    \begin{minipage}[t]{0.5\textwidth}
    \centering\includegraphics{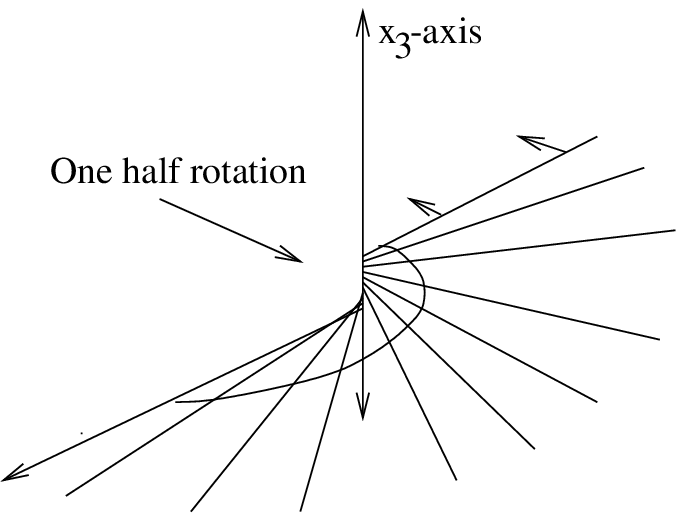}
    \caption{Multi-valued graphs.  The helicoid is obtained
    by gluing together two $\infty$-valued graphs along a line.}\label{f:f1}
    \end{minipage}%
    \begin{minipage}[t]{0.5\textwidth}
    \centering\input{pl8a.pstex_t}
    \caption{The separation $w$ grows/decays in $\rho$ at most sublinearly
for a multi-valued minimal graph; see \eqr{e:sublinerly}.}\label{f:f2}
    \end{minipage}
\end{figure}

\vskip2mm
\noindent
{\bf{Example 2}}:
(Helicoid; see fig. \ref{f:f1}).  The helicoid is the minimal surface 
in $\RR^3$ given by the parametrization 
\begin{equation}  \label{e:helicoid}
(s\cos t,s\sin t,t)\, ,\,\,\,\,\,\text{ where }s,\,t\in \RR\, .
\end{equation}

\vskip2mm
To be able to give a 
precise meaning to the statement 
that the helicoid is a double spiral staircase we will need
the notion of a multi-valued graph, each
staircase will be a multi-valued graph, see fig. \ref{f:f1}.  
Intuitively, an (embedded) multi-valued graph is a 
surface such that over each point of the annulus, 
the surface consists of $N$ graphs. To make this notion 
precise, let $D_r$ be the disk in the plane 
centered at the origin and of radius 
$r$ and let $\cP$ be the universal cover of the
punctured plane $\CC\setminus \{0\}$ with global polar coordinates
$(\rho, \theta)$ so $\rho>0$ and $\theta\in \RR$.  An {\it $N$-valued
graph} on the annulus $D_s\setminus D_r$ is a
single valued graph of a function $u$
over $\{(\rho,\theta)\,|\,r< \rho\leq s\,
,\, |\theta|\leq N\,\pi\}$.  For working purposes, we generally think of the 
intuitive picture of a multi-sheeted surface in 
$\RR^3$, and we identify the single-valued graph over 
the universal cover with its multi-valued image in $\RR^3$.  

The multi-valued graphs that
we will consider will all be
embedded, which corresponds to a nonvanishing separation
between the sheets (or the floors).  Here the {\it separation} is the 
function (see
fig. \ref{f:f2})
\begin{equation}
w(\rho,\theta)=u(\rho,\theta+2\pi)-u(\rho,\theta)\, .
\end{equation}
If $\Sigma$ is the helicoid, then
$\Sigma\setminus \{x_3-\text{axis}\}=\Sigma_1\cup \Sigma_2$,
 where $\Sigma_1$, $\Sigma_2$ are $\infty$-valued graphs on 
$\CC\setminus \{0\}$.
$\Sigma_1$ is the graph of the function $u_1(\rho,\theta)=\theta$
and $\Sigma_2$ is the
graph of the function $u_2(\rho,\theta)=\theta+\pi$.  ($\Sigma_1$ is the 
subset where $s>0$ in \eqr{e:helicoid} and $\Sigma_2$ the subset 
where $s<0$.)  In either
case the separation
$w=2\,\pi$.  A {\it multi-valued minimal graph} is a multi-valued graph of a
function $u$ satisfying the minimal surface equation.

Note that for an embedded multi-valued graph, the sign of $w$ determines
whether the multi-valued graph spirals in a 
left-handed or right-handed manner, in other words, whether upwards motion 
corresponds to turning in a clockwise direction or in a counterclockwise 
direction.  For DNA,
although both spirals occur, the right-handed spiral is far more common 
because of certain details of the chemical structure; see \cite{CaDr}.

As we will see, a fundamental theorem about  embedded minimal disks is that
such a disk is either a minimal graph or can be approximated by a
piece of a rescaled helicoid depending on whether the curvature is small
or not; see Theorem \ref{t:t0.1} below.  To avoid tedious dependence of 
various quantities we state this, our main result, not for a single 
embedded minimal disk with sufficiently large curvature at a given 
point but instead for a sequence of such disks where the curvatures 
are blowing up.  Theorem \ref{t:t0.1} says that a sequence of embedded 
minimal disks 
mimics the following behavior of a sequence of rescaled helicoids. 
 
\parbox{6in}{Consider the sequence 
$\Sigma_i = a_i \, \Sigma$ of 
rescaled helicoids where $a_i \to 0$. (That is, rescale $\RR^3$ by $a_i$, so 
points that used to be distance $d$ apart will in the rescaled $\RR^3$ be 
distance $a_i\,d$ apart.)  The curvatures of this 
sequence of rescaled helicoids 
are blowing up along the vertical axis. The sequence
converges (away from the vertical axis) to a foliation by flat
parallel planes. The singular set $\cS$ (the axis) then consists
of removable singularities.}

\vskip2mm
Throughout let $x_1 , x_2 , x_3$ be the standard coordinates on $\RR^3$.
For $y \in \Sigma \subset \RR^3$ and $s > 0$, the
extrinsic and intrinsic balls are
$B_s(y)$, $\cB_s(y)$.  That is, $B_s(y)=\{x\in \RR^3||x-y|<s\}$ and
$\cB_s(y)=\{x\in \Sigma | \dist_{\Sigma}(x,y)<s\}$.  
$\K_{\Sigma}=\kappa_1\,\kappa_2$ is the Gaussian curvature of
$\Sigma\subset \RR^3$, so when $\Sigma$ is minimal
(i.e., $\kappa_1=-\kappa_2$),
then $|A|^2=-2\,\K_{\Sigma}$.

\vskip2mm
See \cite{CM1}, \cite{O}, \cite{S} (and the forthcoming book \cite{CM3})
for background and basic properties of minimal surfaces and
\cite{CM2} for a more detailed survey of the results described here 
and references.  See also \cite{C} for an abreviated version of this 
paper intended for a general non-mathematical audience.  The article \cite{A} 
discusses in a simple nontechnical way the shape of various things that are 
of ``minimal'' type.  These 
shapes include soap films and soap bubbles, metal alloys, 
radiolarian skeletons, and embryonic tissues and cells.  The reader 
interested in some of the history of the field of 
minimal surfaces may consult \cite{DHKW},  
\cite{N}, and \cite{T}.  

\vskip2mm
\noindent
{\bf{The limit foliation and the singular curve}}
\vskip2mm

\begin{figure}[htbp]
    \setlength{\captionindent}{20pt}
    \begin{minipage}[t]{0.5\textwidth}
    \centering\input{unot3a.pstex_t}
    \caption{Proving Theorem \ref{t:t0.1}.  A. Finding a small 
$N$-valued graph in $\Sigma$.
    B. Extending it in $\Sigma$ to a large
    $N$-valued graph.  C. Extending 
    the number of sheets.}\label{f:f3}
    \end{minipage}\begin{minipage}[t]{0.5\textwidth}
    \centering\input{unot3.pstex_t}
    \caption{Theorem \ref{t:t0.1} - the singular set and the
    two multi-valued graphs.}\label{f:f4}
    \end{minipage}%
\end{figure}

In the next few sections we will discuss how to show that every embedded 
minimal disk is either a graph of 
a function or part of a double spiral staircase; Theorem \ref{t:t0.1} below 
gives precise meaning to this 
statement.  In particular, we will in the next few sections discuss 
the following; 
see fig. \ref{f:f3}:
\vskip1.5mm
\noindent
A.  Fix an integer $N$ (the ``large'' of the curvature in what follows 
will depend on $N$).  If an embedded minimal disk $\Sigma$ is not a 
graph (or equivalently if
the curvature is large at some point), then it contains an
$N$-valued minimal graph which initially  
is shown to exist on the
scale of $1/\max |A|$.  That is, the $N$-valued graph is initially 
shown to be defined on 
an annulus with both inner and outer radius inversely proportional to 
$\max |A|$.
\vskip1.5mm
\noindent
B.  Such a potentially small $N$-valued graph
sitting inside $\Sigma$ can then be seen to extend as an
$N$-valued graph inside $\Sigma$ almost all the way to the
boundary.  That is, the small $N$-valued graph can be extended 
to an $N$-valued graph defined on an annulus where the outer radius of 
the annulus is proportional to $R$.  Here $R$ is the radius of the ball in 
$\RR^3$ that the boundary of $\Sigma$ is contained in. 
\vskip1.5mm
\noindent
C.  The $N$-valued graph not only extends
horizontally (i.e., tangent to the initial sheets) but also
vertically (i.e., transversally to the sheets).  That is, once
there are $N$ sheets there are many more and, in fact, the disk $\Sigma$ 
consists of two multi-valued graphs glued together along an axis.
\vskip1.5mm

These three items, A., B., and C. will be used to demonstrate 
the following theorem, which is the main result:

\begin{Thm} \label{t:t0.1}
(See fig. \ref{f:f4}).
Let $\Sigma_i \subset B_{R_i}=B_{R_i}(0)\subset \RR^3$ be a
sequence of embedded minimal
disks with $\partial \Sigma_i\subset \partial B_{R_i}$
where $R_i\to \infty$. If $\sup_{B_1\cap \Sigma_i}|A|^2\to \infty$, then
there exists a subsequence, $\Sigma_j$,
and 
a Lipschitz curve $\cS:\RR\to \RR^3$ such that after a rotation of $\RR^3$:\\
\underline{1.} $x_3(\cS(t))=t$.  (That is, $\cS$ is a graph 
over the $x_3$-axis.)\\
\underline{2.}  Each $\Sigma_j$ consists of exactly two multi-valued 
graphs away
from $\cS$ (which spiral together).\\
\underline{3.} For each $\alpha>0$, $\Sigma_j\setminus \cS$ converges
in the $C^{\alpha}$-topology to the foliation,
$\cF=\{x_3=t\}_t$, of $\RR^3$ by flat parallel planes.\\ 
\underline{4.}  $\sup_{B_{r}(\cS (t))\cap \Sigma_j}|A|^2\to\infty$ for 
all $r>0$, $t\in \RR$.  (The curvature blows up along $\cS$.)
\end{Thm}

In \underline{2.}, \underline{3.} 
the statement that $\Sigma_j\setminus \cS$ are multi-valued graphs and
converge to $\cF$ means that for each compact subset
$K\subset \RR^3\setminus \cS$
and $j$ sufficiently large $K\cap \Sigma_j$ consists of multi-valued
graphs over (part of) $\{x_3=0\}$
and $K\cap \Sigma_j\to K\cap \cF$.

\vskip2mm
As will be clear in the following sections, A., B., and C. alone are 
not enough to prove Theorem \ref{t:t0.1}.  For instance, \underline{1.} 
does not follow from A., B., and C. but needs a more precise statement than C. 
of where the new sheets form above and below a given multi-valued graph.   
This is done using the ``one-sided curvature estimate''.  

\vskip2mm
Here is a summary of the rest of the paper:

First we discuss two key results that are used in the proof of 
Theorem \ref{t:t0.1}.  These are the existence of
multi-valued graphs, i.e., A. and B., and the important 
one-sided curvature estimate.   Following that we discuss some 
bounds for the separation of multi-valued minimal
graphs.  These bounds are used in both B. and C. above and we discuss what
they are used for in C.  After that we explain how the one-sided 
curvature estimate is
used to show that the singular set, $\cS$, is a Lipschitz curve.  The two last
sections before our concluding remarks 
contain further discussion on the existence of multi-valued graphs and
on the proof of the one-sided curvature estimate.

\vskip2mm
\noindent
{\bf{Two key ingredients in the proof of Theorem \ref{t:t0.1} - 
Existence of multi-valued graphs and the one-sided curvature
estimate}}
\vskip2mm

We now come to the two key results about embedded minimal disks.  The first
says that if the curvature of such a disk $\Sigma$
is large at some point $x\in \Sigma$, then near $x$
 a multi-valued graph forms (in $\Sigma$) and this extends
(in $\Sigma$) almost all the way to the boundary.  
Moreover, the inner radius, $r_x$, of the annulus where the
multi-valued graphs is defined is inversely proportional to $|A|(x)$ 
and the initial separation between the sheets is bounded by a constant times
the inner radius, i.e., $|w(r_x,\theta)|\leq C\,r_x$.

An important ingredient in the proof of Theorem \ref{t:t0.1} is
that, just like the helicoid, general embedded minimal disks with
large curvature at some interior point can be built out of
$N$-valued graphs.  In other words, any
embedded minimal disk can be divided into pieces each of which is
an $N$-valued graph.  Thus the disk itself should be thought of as
being obtained by stacking these pieces (graphs) on top of
each other.  

\begin{figure}[htbp]
    \setlength{\captionindent}{20pt}
    \begin{minipage}[t]{0.5\textwidth}
    \centering\input{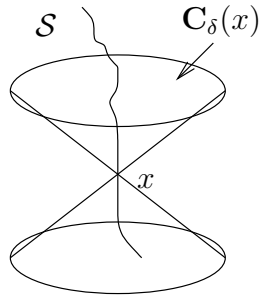}
    \caption{It follows from the one-sided curvature estimate that the singular
    set has the cone property and hence is a Lipschitz curve.}\label{f:f9}
    \end{minipage}
\end{figure}

\begin{figure}[htbp]
    \setlength{\captionindent}{20pt}
\begin{minipage}[t]{0.5\textwidth}
    \centering\input{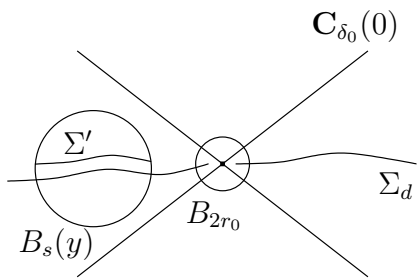}
    \caption{Corollary \ref{c:conecor}:  
     With $\Sigma_d$ playing the role of the plane $x_3=0$,
     by the one-sided estimate, $\Sigma$ consists of 
multi-valued graphs away from a cone.}\label{f:f10}
    \end{minipage}%
\end{figure}

The second key result (Theorem \ref{t:t2}) is a curvature estimate for
embedded minimal disks in a half-space.  As a corollary of this theorem,  
we get that 
the set of points in an embedded minimal disk where the curvature is 
large lies within 
a cone and thus the multi-valued 
graphs, whose existence were discussed above, will all start off within 
this cone; 
see fig. \ref{f:f9} and fig. \ref{f:f10}. 

The curvature 
estimate for disks in a half-space is the following:

\begin{figure}[htbp]
    \setlength{\captionindent}{20pt}
\begin{minipage}[t]{0.5\textwidth}
    \centering\input{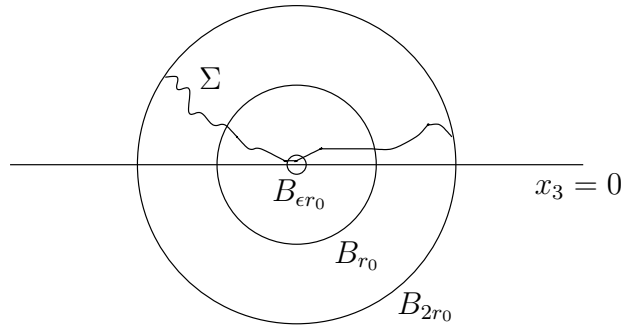}
    \caption{The one-sided curvature estimate for an embedded minimal disk
    $\Sigma$ in a half-space with $\partial \Sigma\subset \partial B_{2r_0}$: 
The components of
    $B_{r_0} \cap \Sigma$ intersecting $B_{\epsilon r_0}$ are graphs.}
\label{f:f6}
    \end{minipage}
\end{figure}

\begin{Thm}  \label{t:t2}
(See fig. \ref{f:f6}).
There exists $\epsilon>0$, such that for all $r_0>0$ if
$\Sigma \subset B_{2r_0} \cap \{x_3>0\}
\subset \RR^3$ is an
embedded minimal
disk with $\partial \Sigma\subset \partial B_{2 r_0}$,
then for all components $\Sigma'$ of
$B_{r_0} \cap \Sigma$ which intersect $B_{\epsilon r_0}$
\begin{equation}  \label{e:graph}
\sup_{x\in\Sigma'} |A_{\Sigma}(x)|^2
\leq r_0^{-2} \, .
\end{equation}
\end{Thm}

Theorem \ref{t:t2} is an interior estimate where the curvature bound, 
\eqr{e:graph}, 
is on the ball $B_{r_0}$ of one half of the radius of the ball $B_{2r_0}$ 
that $\Sigma$ is contained in.  
This is just like a gradient estimate for a harmonic function where the 
gradient bound is on one half of 
the ball where the function is defined.  

Using the minimal surface equation and the fact that $\Sigma'$ has points
close to a plane, it is not hard to see that, for $\epsilon>0$
sufficiently small, \eqr{e:graph} is equivalent to the statement 
 that $\Sigma'$
is a graph over the plane $\{x_3=0\}$.

\begin{figure}[htbp]
    \setlength{\captionindent}{20pt}
    \begin{minipage}[t]{0.5\textwidth}
    \centering\input{unot6.pstex_t}
    \caption{The catenoid given by revolving $x_1= \cosh x_3$
around the $x_3$-axis.}\label{f:f7}
    \end{minipage}\begin{minipage}[t]{0.5\textwidth}
    \centering\input{unot7.pstex_t}
    \caption{Rescaling the catenoid shows that the property of being simply 
    connected (and embedded) is
    needed in the one-sided curvature estimate.}\label{f:f8}
    \end{minipage}%
\end{figure}

We will often refer to Theorem \ref{t:t2}
as {\it the one-sided curvature estimate} (since $\Sigma$ is assumed
to lie on one side of a plane).
Note that the assumption in Theorem \ref{t:t2}
that $\Sigma$ is simply connected (i.e., that $\Sigma$ is a disk) is crucial
as can be seen from the example of a rescaled catenoid. The catenoid,
see fig. \ref{f:f7},
is the minimal surface in $\RR^3$ given by
$(\cosh s\, \cos t,\cosh s\, \sin t,s)$
where $s,t\in\RR$.
Rescaled catenoids converge (with multiplicity two) to
the flat plane; see fig. \ref{f:f8}.  Likewise, by considering
the universal cover of the catenoid, one sees that Theorem \ref{t:t2} requires 
the disk to be embedded,
and not just immersed.

\begin{Def}   \label{d:cone}
(Cones; see fig. \ref{f:f9}).  
If $\delta>0$ and $x\in \RR^3$, then we denote by
$\cone_{\delta}(x)$ the (convex) cone with vertex $x$, cone angle $(\pi/2 -
\arctan \delta)$, and axis parallel to the $x_3$-axis.    That is,
(see fig. \ref{f:f9})
\begin{equation}
\cone_{\delta}(x)=\{x\in \RR^3\,|\,x_3^2 \geq
\delta^2\,(x_1^2+x_2^2) \} +x\, .
\end{equation}
\end{Def}

In the proof of Theorem \ref{t:t0.1},
the following (direct) consequence of Theorem \ref{t:t2} 
(with
$\Sigma_d$ playing the role of the plane $\{x_3=0\}$, see fig. 
\ref{f:f10}) is needed 
(Paraphrased, this corollary says that if an embedded minimal disk contains a
$2$-valued graph, then the disk consists of multi-valued graphs away from a
cone with axis orthogonal to the $2$-valued graph.  In Corollary 
\ref{c:conecor} the ``d'' in $\Sigma_d$ stands for double-valued.):

\begin{Cor}   \label{c:conecor}
(See fig. \ref{f:f10}).
There exists $\delta_0>0$ so for all $r_0$, $R>0$ with $r_0<R$ 
if $\Sigma\subset B_{2R}$,
$\partial \Sigma\subset \partial B_{2R}$ is an embedded minimal
disk containing a $2$-valued graph $\Sigma_d \subset
\RR^3\setminus \CC_{\delta_0}(0)$ over the annulus 
$D_{R}\setminus D_{r_0}$ with
gradient $\leq \delta_0$, then each component of $B_{R/2}\cap
\Sigma\setminus (\cone_{\delta_0}(0)\cup B_{2 r_0})$ is a
multi-valued graph.
\end{Cor}

Fig. \ref{f:f10} illustrates how this corollary follows from Theorem
\ref{t:t2}.  In this picture, $B_s(y)$ is a ball away from $0$
and $\Sigma'$ is a component of $B_s(y)\cap \Sigma$ disjoint from
$\Sigma_d$.  It follows easily from the maximum principle that
$\Sigma'$ is topologically a disk.  Since $\Sigma'$ is assumed to
contain points near $\Sigma_d$, then we can let a component of
$B_s(y)\cap \Sigma_d$ play the role of the plane $\{x_3=0\}$ in
Theorem \ref{t:t2} and the corollary follows.

\vskip2mm
Using Theorems \ref{t:t0.1}, \ref{t:t2},
W. Meeks and H. Rosenberg proved in ``The uniqueness of the helicoid and 
the asymptotic geometry of properly
embedded minimal surfaces with finite topology'' 
that the plane and helicoid are the only
complete properly embedded
simply-connected minimal surfaces in $\RR^3$.   
Catalan had proven in 1842 that any complete ruled minimal surface is 
either a plane or a helicoid.  
A surface is said to be {\it ruled} if it has the 
parametrization   
\begin{equation}
X(s,t)=\beta (t)+ s\,\delta (t)\, ,\,\,\,\,\,\text{ where }s,\,t\in \RR\, ,
\end{equation}
and $\beta$, $\delta$ are curves in $\RR^3$.  The curve $\beta (t)$ 
is called the {\it directrix} of the surface, and a line having 
$\delta (t)$ as direction vector is called a {\it ruling}.  For the 
helicoid in \eqr{e:helicoid}, the $x_3$-axis is a directrix, 
and for each fixed $t$ the line $s\to (s\,\cos t,s\,\sin t,t)$ is a ruling.  

\vskip2mm
\noindent
{\bf{Towards removablility of singularities - Analysis of 
multi-valued minimal graphs}}
\vskip2mm

Even given the decomposition into multi-valued graphs mentioned in the 
beginning of the 
previous section, to prove Theorem
\ref{t:t0.1}, one still needs to analyze how the various $N$-valued
pieces fit together.  In particular, we need 
Theorem \ref{t:t2} and Corollary \ref{c:conecor} to show that an embedded 
minimal disk that is not a graph 
cannot be contained in a half-space and thus the subset of points with 
large curvature lies within a cone.  This is
still not
enough to imply Theorem \ref{t:t0.1}. One also needs to show
that part of any embedded minimal disk cannot accumulate in a
half-space.  This is what we call {\it properness} below; see fig. 
\ref{f:f5} and
\eqr{e:exofu} that
gives an example of an $\infty$-valued graph whose image lies in a slab in 
$\RR^3$.  The
property we call properness is the assertion that no limit of 
embedded minimal disks 
can contain
such a (nonproper) multi-valued graph.  

In this section, we will discuss bounds for the separation of embedded 
multi-valued graphs and their applications to properness and to the proofs of 
Theorems \ref{t:t0.1}, \ref{t:t2}.  Two types of bounds for the 
growth/decay (as
$\rho\to \infty$) of the separation will be needed:
\vskip1.5mm
\noindent
a.  The
weaker sublinear bounds, i.e., there exists $0<\alpha<1$ such that 
for fixed $\rho_0$, we have the
bounds 
\begin{equation}   \label{e:sublinerly}
(\rho/\rho_0)^{-\alpha}\,|w(\rho_0,\theta)|\leq |w(\rho,\theta)|\leq
(\rho/\rho_0)^{\alpha}\,|w(\rho_0,\theta)| 
\,\,\,\,\,\text{ as }\rho\to \infty\, .
\end{equation}
These bounds hold for $N$-valued graphs (where $N$ is
some fixed large number).   By letting $N$ be large, 
$\alpha$ can be chosen small.
\vskip1.5mm
\noindent
b. The stronger logarithmic bounds, i.e., there exist constants 
$c_1$ and $c_2$ such that for fixed $\rho_0$, we have the
bounds  
\begin{equation}  \label{e:logarithmicbound}
\frac{c_1}{\log (\rho/\rho_0)}\,|w(\rho_0,\theta)|\leq |w(\rho,\theta)|\leq
c_2\,\log (\rho/\rho_0)\,|w(\rho_0,\theta)|\,\,\,\,\,  
\text{ as }\rho\to \infty\, .
\end{equation}
These bounds 
will require a growing number of sheets (growing as
$\rho\to\infty$) and will be used only to show properness; cf. fig. \ref{f:f5}.

\begin{figure}[htbp]
    \setlength{\captionindent}{20pt}
    \begin{minipage}[t]{0.5\textwidth}
    \centering\input{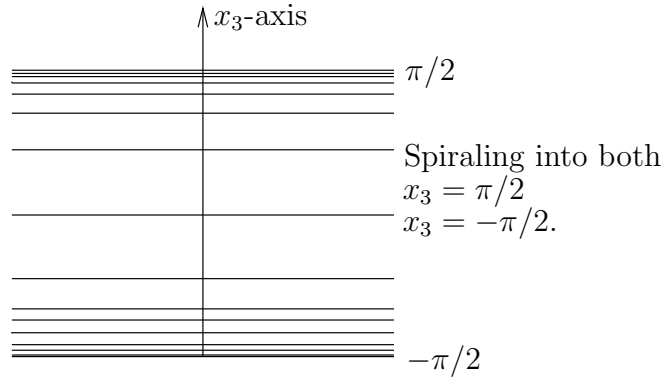}
    \caption{To show properness one needs to rule out
    that one of the multi-valued graphs  can
    contain a nonproper graph  like  $\arctan (\theta / \log \rho)$,
    where $(\rho , \theta)$ are polar coordinates.  The graph of 
    $\arctan (\theta / \log \rho)$ 
    is illustrated above.}\label{f:f5}
    \end{minipage}
\end{figure}

\vskip1.5mm
Here are couple of things that the sublinear bounds are used for.  
First, as a consequence of the existence of multi-valued graphs 
discussed in the previous section,  
one easily
gets that if $|A|^2$ is blowing up near $0$ for a sequence of
embedded minimal
disks $\Sigma_i$,
then there is a sequence of $2$-valued graphs
$\Sigma_{i,d}\subset \Sigma_i$.  Here the $2$-valued graphs start
off defined on a smaller and smaller scale (the inner radius of the 
annulus where each 
multi-valued graph is defined is inversely proportional to $|A|$).
Consequently, by the sublinear separation growth, such $2$-valued
graphs collapse:  Namely, if
$\Sigma_{i,d}$ is a $2$-valued graph over $D_R\setminus D_{r_i}$, then
$|w_i(\rho,\theta)|\leq (\rho/r_i)^{\alpha}\,|w_i(r_i,\theta)|\leq
C\,\rho^{\alpha}\,r_i^{1-\alpha}$ for some $\alpha<1$ and some 
constant $C$.  (In fact, by making $N$ large, $\alpha$ can be chosen small.) 
Letting $r_i\to 0$ shows that $|w_i(\rho,\theta)|\to 0$ for $\rho$, 
$\theta$ fixed.
Thus as $i\to\infty$ the upper sheet collapses onto the lower 
and, hence, a subsequence converges to a smooth
minimal graph through $0$.  (Here $0$ is a removable singularity
for the limit.)  Moreover, if the sequence of such disks
is as in Theorem
\ref{t:t0.1}, i.e., if $R_i\to \infty$, then the minimal graph in the
limit is entire and hence, by Bernstein's theorem, is a plane.

The sublinear bounds are also used in the proof of Theorem \ref{t:t2} 
which
in turn - through its corollary, Corollary \ref{c:conecor} -  
is used to show that any multi-valued graph contained 
in an embedded minimal
disk can be extended, inside the disk, to a multi-valued graph 
with a rapidly growing number of sheets and thus we get the
better logarithmic bounds for the separation.  Namely, by Corollary 
\ref{c:conecor}, 
outside a cone such a multi-valued graph extends as a multi-valued graph.  
An application of a Harnack inequality then shows that the number of 
sheets that it takes to leave 
the complement of the cone where the disk is graphical must grow in $\rho$ 
sufficiently fast so that \eqr{e:logarithmicbound} follows.  (Recall that 
for the bounds \eqr{e:logarithmicbound} to hold requires that the number 
of sheets grows sufficiently fast.)   

By b. when the number of sheets grows sufficiently fast, 
the fastest possible decay for
$w(\rho,0)/w(1,0)$ is $c_1/\log \rho$.  
This lower bound for the decay of the separation 
is sharp.  It is 
achieved for the $\infty$-valued graph of the harmonic function (graphs of
multi-valued 
harmonic functions are good models for multi-valued minimal graphs)
\begin{equation}  \label{e:exofu}
u(\rho,\theta)=\arctan  \frac{\theta}{\log \rho} \, .
\end{equation}
Note that the graph of $u$ is embedded and lies in a slab in $\RR^3$, 
i.e., $|u|\leq \pi/2$,
and hence in particular is not proper.  On the top it spirals into the
plane $\{x_3=\pi/2\}$ and on the bottom into $\{x_3=-\pi/2\}$, yet it never
reaches either of these planes; see fig. \ref{f:f5}.  

The next proposition rules out not only this as a possible limit of
(one half of)  embedded
minimal disks, but, more generally, any
$\infty$-valued minimal graph in
a half-space.

\begin{Pro}  \label{c:main}
Multi-valued graphs contained in embedded minimal disks 
are proper - they do not
accumulate in finite height.
\end{Pro}

Proposition \ref{c:main} relies in part on the logarithmic bound, 
\eqr{e:logarithmicbound}, 
for the separation.  

\vskip2mm
\noindent
{\bf{Regularity of the singular set and Theorem \ref{t:t0.1}}}
\vskip2mm

In this section we will indicate how to define the singular set $\cS$ 
in Theorem \ref{t:t0.1} and show regularity of $\cS$.  

By a very general standard compactness argument, it follows (after possibly
going to a subsequence)
that for a sequence of smooth surfaces
there is a well defined notion of points where the second
fundamental form of the sequence blows up.  That is,
let $\Sigma_i\subset B_{R_i}$, $\partial \Sigma_i\subset \partial B_{R_i}$,
and $R_i\to \infty$ be a sequence of (smooth) compact surfaces.
After passing to a subsequence, $\Sigma_j$, we
may assume that for each $x\in \RR^3$ either (a) or (b) holds:\\
(a) $\sup_{B_{r}(x)\cap
\Sigma_j}|A|^2\to \infty$ for all $r>0$,\\
(b) $\sup_j\sup_{B_r(x)\cap
\Sigma_j}|A|^2<\infty$ for some $r>0$.

\begin{Def}
(Cone property).  
Fix $\delta>0$. We will say that a subset $\cS\subset \RR^3$ has the
{\it cone property} (or the $\delta$-cone property) if $\cS$ is closed and
nonempty and:\\
(i) If $z\in \cS$, then $\cS\subset
\cone_{\delta}(z)$; see Definition \ref{d:cone} for $\CC_{\delta}(z)$.\\
(ii) If $t\in x_3(\cS)$ and
$\epsilon>0$, then $\cS\cap \{t<x_3<t+\epsilon\}\ne \emptyset$ and
$\cS\cap \{t-\epsilon<x_3<t\}\ne \emptyset$.
\end{Def}

\vskip2mm
Note that (ii) just says that each point in $\cS$ is the limit of points
coming from above and below.

When $\Sigma_i\subset B_{R_i}\subset \RR^3$  is a sequence of embedded
minimal disks with $\partial \Sigma_i\subset \partial B_{R_i}$,
$R_i\to \infty$, $\Sigma_j$ is the subsequence as above,
and $\cS$ is the set of points where the curvatures of $\Sigma_j$
blow up (i.e., where (a) above holds), then we will see below that
$\cS$ has the
cone property (after a rotation of $\RR^3$).   Hence (by the next lemma),
$\cS$ is a Lipschitz curve which
is a graph over the $x_3$-axis.  Note that when $\Sigma_i$
is a sequence of rescaled helicoids, then $\cS$ is the $x_3$-axis.

\begin{Lem}  \label{l:regsing}
(See fig. \ref{f:f9}).
If $\cS\subset \RR^3$ has the $\delta$-cone property, then
$\cS\cap
\{x_3=t\}$ consists of exactly one point $\cS_t$ for all
$t\in\RR$, and $t\to \cS_t$ is a Lipschitz parameterization of
$\cS$.
\end{Lem}

With Lemma \ref{l:regsing} in hand, we can proceed with the proof of 
Theorem \ref{t:t0.1}.  So suppose that $\Sigma_i$ is as in Theorem 
\ref{t:t0.1} and
$\Sigma_j$, $\cS$ are as above, then $\cS$ is closed by definition
and nonempty by the assumption of Theorem \ref{t:t0.1}.  Centered at
any $x\in \cS$ we can, by the existence of multi-valued graphs
near points where the curvatures blow up, the sublinear separation growth,
and Bernstein's theorem, find a sequence of $2$-valued graphs 
$\Sigma_{d,j}\subset \Sigma_j$ which
converges to a plane through $x$; see the discussion
preceding Theorem \ref{t:t2}.  (This is after possibly passing to a
subsequence of the $\Sigma_j$'s.)  Thus (i) above
holds by Corollary \ref{c:conecor}.  Therefore to
see that $\cS$ has the cone property all we need to see is that
(ii) holds.  The proof of this relies on Proposition \ref{c:main}.
Once the cone property of $\cS$ is shown, it follows from
Lemma  \ref{l:regsing} that $\cS$ is a Lipschitz curve and by
Corollary \ref{c:conecor}, away from $\cS$, each $\Sigma_j$ consists of
multi-valued graphs.  It is not hard to see that there are at
least two such graphs and a barrier argument shows that
there are not more.

\vskip2mm
\noindent
{\bf{Blow up points and the existence of multi-valued graphs}}
\vskip2mm

To describe the existence of
multi-valued graphs in embedded minimal disks,
we will need the notion of a blow up point.

Let $\Sigma$ be a smooth (minimal or not) embedded 
(compact) surface in a ball $B_{r_0}(x)$ in $\RR^3$, passing through $x$ 
the center of the ball, and with boundary contained in the 
boundary of the ball.  Here $B_{r_0}(x)$ is the
extrinsic ball of radius $r_0$, but could as well have been an
intrinsic ball $\cB_{r_0}(x)$ 
in which case the notion of a blow up point below
would have to be appropriately changed. Suppose that $|A|^2(x)\geq
4 \, C^2\,r_0^{-2}$ for some constant $C>0$.  We claim that there is
$y\in B_{r_0}(x)\cap \Sigma$ and $s>0$ such that $B_s(y)\subset
B_{r_0}(x)$ and
\begin{equation}  \label{e:ti}
\sup_{B_s(y)\cap \Sigma}|A|^2\leq 4\,C^2\,s^{-2}= 4\,|A|^2(y)\, .
\end{equation}
That is, the curvature at $y$ is large (this just means that $C$ should
be thought of as a large constant equal to $s\,|A|(y)$) and is almost
(up to the constant $4$) the maximum on
the ball $B_s(y)$.  We will say that the pair $(y,s)$ is a 
{\it blow up pair} and the point $y$ is a {\it blow up point}.   
Later $s$ will be replaced by $8s$ and eventually by a constant times $s$ 
and sometimes the extrinsic ball will be replaced by an intrinsic ball, but 
we will still refer to the pair $(y,s)$ as a blow up pair.   
That there exists such a point $y$ is easy to see;
on $B_{r_0}(x)\cap \Sigma$ set $F(z)=(r_0-r(z))^2\,|A|^2(z)$ where
$r(z)=|z-x|$.  Then 
\begin{equation}
F(x)\geq 4\,C^2\, , F\geq 0\, , \text{ and }
\left. F \right|_{\partial B_{r_0}(x)\cap \Sigma}=0\, .  
\end{equation}
Let $y$ be where the maximum
of $F$ is achieved and set $s= C/|A|(y)$.  Note that $s$ is at most 
one-half of the distance from $y$ to the boundary of the ball $B_{r_0}(x)$.  
One easily checks that
$y$, $s$ have the required properties.  Namely, clearly $|A|^2(y)= C^2\,s^{-2}$
and since $y$ is where the maximum of $F$ is achieved, 
\begin{equation}
|A|^2(z)\leq
\left(\frac{r_0-r(y)}{r_0-r(z)}\right)^2\,|A|^2(y)\, .  
\end{equation}
Since $F(x)\geq 4\,C^2$ it follows from
the choice of $s$ that $|r_0-r(y)|\leq 2\,|r_0-r(z)|$ for 
$z\in B_s(y)\cap \Sigma$.  Hence,
$|A|^2(z)\leq 4\,|A|^2(y)$.  Together this gives \eqr{e:ti}.  

The existence of multi-valued graphs
is shown 
by combining a blow up result with an extension result.
This blow up result says that if an embedded
minimal disk in a ball has large curvature at a point, then it
contains a small (in fact on the scale of $1/\max |A|$) 
almost flat $N$-valued graph nearby; this is A. in
fig. \ref{f:f3}.
The extension result allows us to extend the (small) $N$-valued graphs
almost out to the boundary of the ``big'' ball $B_R$; this is B. in
fig. \ref{f:f3}.  
In fact, the blow up result shows that if $(y,s)$ is a blow up pair
with point $y$ and radius $s>0$ satisfying \eqr{e:ti}, then the corresponding
$N$-valued function is defined on an annulus whose inner radius is $s$ and so
the initial separation is proportional to $s$.  
That is, for positive constants 
$C_1$, $C_2$
\begin{equation}  \label{e:proport}
C_1\, s\leq |w(s,\theta)|\leq C_2\, s\, .
\end{equation}
Equation \eqr{e:proport} will be used implicitly in the next section.

\begin{figure}[htbp]
    \setlength{\captionindent}{20pt}
    \begin{minipage}[t]{0.5\textwidth}
    \centering\input{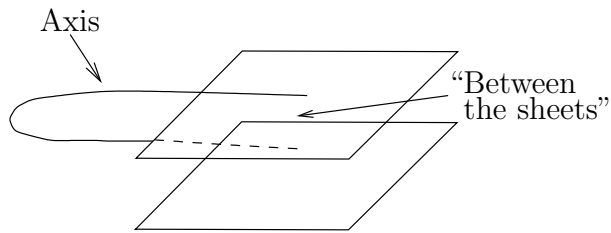}
    \caption{The curvature estimate \newline ``between the sheets.''}
    \label{f:f11}
    \end{minipage}
\end{figure}

The extension result is significantly more subtle than
the local existence of multi-valued graphs.  
The key for being able to extend is a curvature
estimate ``between the sheets'' for embedded minimal disks; 
see fig. \ref{f:f11}.
We think of an axis for such a disk $\Sigma$ as a point or
curve away from which the
surface locally (in an extrinsic ball) has more than one
component.  With this weak notion of an axis, the estimate between the sheets 
is that if
one component of $\Sigma$ is sandwiched between two
others that connect to an axis, then there is a fixed bound for (the norm
of) 
the curvature of the one that
is sandwiched.  
The example to keep in
mind is a helicoid and the components are
``consecutive sheets'' away from the axis.  Once the estimate between
the sheets is established, then it is applied to the ``middle'' sheet(s)
of an $N$-valued graph to show that even as we go far out to the ``outer''
boundary of the $N$-valued graph the
curvature has a fixed bound.  Using this a priori bound and
additional arguments one 
gets better bounds and eventually (with more work) argue that the
sheets must remain almost flat and thus the 
$N$-valued graph will remain an $N$-valued graph.

\vskip2mm
\noindent
{\bf{The proof of the one-sided curvature estimate}}
\vskip2mm

Using a blow up argument and the minimal surface equation,
one can show curvature estimates for minimal
surfaces which on all
sufficiently small scales lie on one side of, but come close to,
a plane. Such an assumption is a scale invariant version of Theorem
\ref{t:t2}.  However, the assumption of Theorem \ref{t:t2} is not scale
invariant and the theorem cannot be proven this way.
The scale invariant condition is very similar
to the classical Reifenberg property.  (A subset of $\RR^n$
has the Reifenberg property if it is close on all scales to a hyper-plane; see
the appendix of \cite{ChC}.  This property goes back to Reifenberg's 
fundamental 1960 paper ``Solution of the Plateau problem for 
$m$-dimensional surfaces of varying topological type''.)
As explained above (in particular Corollary \ref{c:conecor}),
the significance of Theorem \ref{t:t2}
is indeed that it only requires closeness on one scale.  On the other
hand, this is what makes it difficult to prove (the lack of scale
invariance is closely related to the lack of a useful monotone
quantity).

Let us give a very rough outline of the proof of the one-sided
curvature estimate; i.e.,
Theorem \ref{t:t2}. Suppose that $\Sigma$
is an embedded minimal disk in the half-space $\{x_3> 0\}$ intersected 
with the ball $B_{2r_0}$ and with boundary in the boundary of the ball 
$B_{2r_0}$. The
curvature estimate is proven by contradiction; so suppose that $\Sigma$ has low
points with large curvature.
Starting at such a point, we decompose
$\Sigma$ into disjoint multi-valued graphs using the existence
of nearby points with large curvature (the existence of such nearby
points is highly nontrivial to establish.  We will use that such
a nearby point of large curvature can be found below any given multi-valued
graph and thus we can choose the ``next'' blow up point to always 
be below the
previous). 
The key point is then to show
(see Proposition \ref{p:lift} below)
that we can in fact find such a decomposition where the ``next''
multi-valued graph starts off a definite amount below where the
previous multi-valued graph started off. In fact, what we show
is that this definite amount is a fixed fraction of the distance
between where the two graphs started off. Iterating this
eventually forces $\Sigma$ to have points where $x_3<0$.
This is the desired contradiction.

\vskip2mm
To show this key proposition
(Proposition \ref{p:lift}), we use two decompositions
and two kinds of blow up points.
The first decomposition uses
the more standard blow up points given as pairs $(y,s)$
where $y\in \Sigma$ and $s>0$ is such that
\begin{equation}   \label{e:defc1ii}
\sup_{\cB_{8s}(y)}|A|^2\leq 4|A|^2(y)=4\,C_1^2\,s^{-2}\, .
\end{equation}
(Here $\cB_{8s}(y)$ is the intrinsic ball of radius $8s$, so in 
particular $\cB_{8s}(y)\subset \Sigma$.)  The point about such a 
pair $(y,s)$ is that, 
since $\Sigma$ is a minimal disk, $\Sigma$ 
contains a multi-valued graph near $y$ starting
off on the scale $s$.  (This is assuming that the curvature at $y$ is
sufficiently large, i.e., $C_1$ is a sufficiently
large constant.)  The second kind of blow up
points are the ones where (except for a technical issue)
$8$ in the radius of the ball centered at $y$ 
is replaced by some really large
constant $C$, i.e.,
\begin{equation}  \label{e:pairs}
\sup_{\cB_{Cs}(y)}|A|^2\leq 4|A|^2(y)=4\,C_1^2\,s^{-2}\, .
\end{equation}

\begin{figure}[htbp]
    \setlength{\captionindent}{20pt}
    \begin{minipage}[t]{0.5\textwidth}
    \centering\input{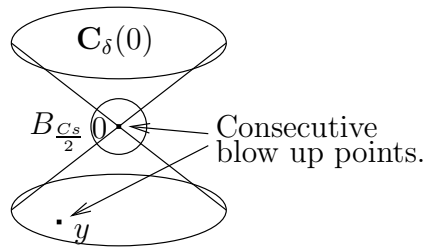}
    \caption{Two consecutive blow up points satisfying \eqr{e:pairs}.}
    \label{f:f12}
    \end{minipage}
\end{figure}

\begin{Pro}   \label{p:lift}
(See fig. \ref{f:f12}).
There exists $\delta>0$ such that if point $0$, radius $s$ 
satisfies \eqr{e:pairs} and $\Sigma_0\subset \Sigma$ is the corresponding
(to $(0,s)$) $2$-valued graph over the annulus $D_{r_0}\setminus D_s$, then
we get point $y$, radius $t$ satisfying \eqr{e:pairs} 
with $y\in \cone_{\delta}(0)
\cap \Sigma \setminus B_{Cs/2}$ and where $y$ is below $\Sigma_0$.
\end{Pro}

\begin{figure}[htbp]
    \setlength{\captionindent}{20pt}
    \begin{minipage}[t]{0.5\textwidth}
    \centering\input{unot11.pstex_t}
    \caption{Between two consecutive blow up points satisfying  \eqr{e:pairs}
there are many blow up points satisfying
 \eqr{e:defc1ii}.}\label{f:f13}
    \end{minipage}%
    \begin{minipage}[t]{0.5\textwidth}
    \centering\input{unot12.pstex_t}
    \caption{Measuring height.  Blow up points and
corresponding multi-valued graphs.}\label{f:f14}
    \end{minipage}
\end{figure}

The point for proving Proposition \ref{p:lift} is 
that we can find blow up points satisfying
\eqr{e:pairs} so that the distance between them is proportional to
the sum of the scales.  Moreover, between consecutive blow up
points satisfying \eqr{e:pairs}, we can find many blow up
points satisfying \eqr{e:defc1ii}; see fig. \ref{f:f13}.  The advantage is
now that if we look between blow up points satisfying
\eqr{e:pairs}, then the height of the multi-valued graph given by
such a pair grows like a small power of the distance whereas the
separation between the sheets in a multi-valued graph given by
\eqr{e:defc1ii} decays (at the worst) like a small power of the distance; see
fig. \ref{f:f14}.  Now, thanks to that the number of blow up points
satisfying \eqr{e:defc1ii} (between two consecutive blow up points
satisfying \eqr{e:pairs}) grows almost linearly, then, even though
the height of the graph coming from the blow up point satisfying
\eqr{e:pairs} could move up (and thus work against us), the sum of
the separations of the graphs coming from the points satisfying
\eqr{e:defc1ii} dominates the other term.  Thus the next blow up
point satisfying \eqr{e:pairs} (which lies below all the other
graphs) is forced to be a definite amount lower than the previous
blow up point satisfying \eqr{e:pairs}.  This gives the proposition.  

Theorem \ref{t:t2} follows from the proposition.  Suppose the theorem
fails; starting at a point of large curvature and iterating the
proposition will eventually give a point in the minimal surface with
$x_3<0$, which is a contradiction.

\vskip2mm
\noindent
{\bf{Concluding remarks and some possible future directions of research}}
\vskip2mm

In this article we have seen that minimal surfaces 
and surfaces of ``minimal'' type occur naturally,  
and we have described why embedded 
minimal disks are double spiral staircases.  
We would hope that similar considerations can be used to 
answer, for other things than minimal disks, the age-old question: 
 
\vskip2mm
\centerline{``What are the possible shapes of various things and why?''.}  
\vskip2mm

A different possible direction is to describe $3$-manifolds.  
Namely, by a result of B. White for a generic metric on a 
closed $3$-manifold the area functional, 
that is, the map that to each closed surface assigns its area, is a 
Morse function.  (Recall that a Morse function is a function that only 
has nondegenerate critical points.)  As we saw earlier, the  
critical points of the area functional are precisely the minimal surfaces;  
thus if one could understand 
all minimal surfaces in a given $3$-manifold, $M$, then one would 
understand all critical points for the area functional.  For a generic metric 
one could 
then hope to use Morse theoretic arguments to understand $M$.  For 
general embedded minimal surfaces of a given fixed genus in closed 
$3$-manifold, the key for understanding them is to understand 
their intersection with a small ball in $M$.  Since locally any fixed 
$3$-manifold is almost euclidean this boils down to understanding minimal 
surfaces in a ball in $\RR^3$ with boundary in the boundary of the ball. The 
key for this is indeed the case where the minimal surfaces are disks, thus 
the key is the results described here.  We will discuss this elsewhere. 

The field of minimal surfaces has undergone enormous development 
since the days of Euler and Lagrange. It has played a key role in the 
development of many other fields in analysis and geometry.  It has had 
times of intense development followed by times of stagnation before new 
fundamental results and techniques have been discovered.  This has happened  
over and over again.  
In closing, we believe that, after nearly a quarter of a millenium,  
the field of minimal surfaces 
is at its very peak and of utmost importance in mathematics and its 
applications.  We hope that this expository article has 
helped convey this.  Although 
as the saying goes ``it is hard to predict -- 
especially about the future''
\footnote{Quote attributed to Danish humorist Storm P. 
(Robert Storm Petersen) from the 1920s.}, 
we believe that more 
magnificent results and techniques 
are to follow.  

\vskip2mm
We are grateful to Christian Berg, Cornelius H. Colding, Chris Croke, 
Bruce Kleiner, Camillo De Lellis, Paul Schlapobersky, David Ussery, and 
David Woldbye  for suggestions and 
comments, and we are particularly grateful to Andrew Lorent for his 
very many very helpful comments.
Finally, we wish to thank the editor of the Notices of the AMS, 
Harold Boas, for his and the reviewer's considerable help 
in making this article more readable.

\end{document}